\documentclass[12pt,leqno]{amsart}
\usepackage{amsmath,amsfonts,amssymb,amsthm}
\usepackage{accents}
\usepackage{caption}
\usepackage{subcaption}
\usepackage{calligra}

\numberwithin{figure}{section}

\theoremstyle{plain}
\newtheorem{theorem}{Theorem}[section]

\theoremstyle{definition}
\newtheorem{defn}{Definition}[section]

\theoremstyle{remark}

\usepackage{mathtools}
\newcommand{\KN}{\mathbin{\bigcirc\mspace{-15mu}\wedge\mspace{3mu}}}
\title[Compact gradient $\rho$-Einstein soliton]{Compact gradient $\rho$-Einstein soliton is isometric to the Euclidean sphere}
\author[A. A. Shaikh, C. K. Mondal, P. Mandal ]{Absos Ali Shaikh$^1$, Chandan Kumar Mondal$^2$, Prosenjit Mandal$^3$ }
\address{\noindent\newline Department of Mathematics,\newline University of
Burdwan, Golapbag,\newline Burdwan-713104,\newline West Bengal, India}
\email{aask2003@yahoo.co.in, aashaikh@math.buruniv.ac.in}
\email{chan.alge@gmail.com}
\email{prosenjitmandal235@gmail.com}
\begin{document}
\begin{abstract}
In this paper we have investigated some aspects of gradient $\rho$-Einstein Ricci soliton in a complete Riemannian manifold. First, we have proved that the compact gradient $\rho$-Einstein soliton is isometric to the Euclidean sphere by showing that the scalar curvature becomes constant. Second, we have showed that in a non-compact gradient $\rho$-Einstein soliton satisfying some integral condition, the scalar curvature vanishes.
\end{abstract}
\noindent\footnotetext{$\mathbf{2020}$\hspace{5pt}Mathematics\; Subject\; Classification: 53C20; 53C21; 53C44. }
\maketitle
\section{Introduction and preliminaries}
A $1$-parameter family of metrics $\{g(t)\}$ on a Riemannian manifold $M$, defined on some time interval $I\subset\mathbb{R}$ is said to satisfy Ricci flow if it satisfies
$$\frac{\partial}{\partial t}g_{ij}=-2R_{ij},$$
where $R_{ij}$ is the Ricci curvature with respect to the metric $g_{ij}$. Hamilton \cite{HA82} proved that for any smooth initial metric $g(0)=g_0$ on a closed manifold, there exists a unique solution $g(t)$, $t\in [0,\epsilon)$, to the Ricci flow equation for some $\epsilon>0$. A solution $g(t)$ of the Ricci flow of the form
$$g(t)=\sigma(t)\varphi(t)^*g(0),$$
where $\sigma:\mathbb{R}\rightarrow\mathbb{R}$ is a positive function and $\varphi(t):M\rightarrow M$ is a 1-parameter family of diffeomorphisms, is called a Ricci soliton. It is known that if the initial metric $g_0$ satisfies the equation
\begin{equation}\label{eq6}
Ric(g_0)+\frac{1}{2}\pounds_Xg_0=\lambda g_0,
\end{equation}
where $\lambda$ is a constant and $X$ is a smooth vector field on $M$, then the manifold $M$ admits Ricci soliton. Therefore, the equation (\ref{eq6}), in general, is known as Ricci soliton. If $X$ is the gradient of some smooth function, then it is called gradient Ricci soliton. For more results of Ricci soliton see \cite{CA10,CK04,FMZ08}. In 1979, Bourguignon \cite{BOU81} introduced the notion of Ricci-Bourguignon flow, where the metrics $g(t)$ is evolving according to the flow equation
$$\frac{\partial}{\partial t}g_{ij}=-2R_{ij}+2\rho Rg_{ij},$$
where $\rho$ is a non-zero scalar constant and $R$ is the scalar curvature of the metric $g(t)$. Following the Ricci soliton, Catino and Mazzier \cite{CM16} gave the definition of gradient $\rho$-Einstein soliton, which is the self-similar solution of Ricci-Bourguignon flow. This soliton is also called gradient Ricci-Bourguignon soliton by some authors.
\begin{defn}\cite{CM16}
Let $(M,g)$ be a Riemannian manifold of dimension $n$, $(n\geq 3)$, and let $\rho\in\mathbb{R}$, $\rho\neq 0$. Then $M$ is called gradient $\rho$-Einstein soliton, denoted by $(M,g,f,\rho)$, if there is a smooth function $f:M\rightarrow\mathbb{R}$ such that
\begin{equation}\label{eq1}
Ric+\nabla^2 f=\lambda g+\rho R g,
\end{equation}
for some constant $\lambda$.
\end{defn}
The soliton is trivial if $\nabla f$ is a parallel vector field. The function $f$ is known as $\rho$-Einstein potential function. If $\lambda>0$ $(\text{resp.} =0,<0 )$, then the gradient $\rho$-Einstein soliton $(M,g,f,\rho)$ is said to be shrinking (resp. steady or expanding) . On the other hand, the $\rho$-Einstein soliton is called gradient Einstein soliton, gradient traceless Ricci soliton or gradient Schouten soliton if $\rho=1/2,1/n$ or $1/2(n-1)$. Later, this notion has been generalized in various directions such as $m$-quasi Einstein manifold \cite{HLX15}, $(m,\rho)$-quasi Einstein manifold \cite{HW13}, Ricci-Bourguignon almost
soliton \cite{DW18}. \\
\indent Catino and Mazzier \cite{CM16} showed that compact gradient Einstein, Schouten or traceless Ricci soliton is trivial. They classified three-dimensional gradient shrinking Schouten soliton and proved that it is isometric to a finite quotient of either $\mathbb{S}^3$ or $\mathbb{R}^3$ or $\mathbb{R}\times\mathbb{S}^2$. Huang \cite{HU17} deduced a sufficient condition for the compact gradient shrinking $\rho$-Einstein soliton to be isometric to a quotient of the round sphere $\mathbb{S}^n$. 
\begin{theorem}\cite{HU17}
Let $(M,g,f,\rho)$ be an $n$-dimensional $(4\leq n\leq 5)$ compact gradient shrinking $\rho$-Einstein soliton with $\rho<0$. If the following condition holds
\begin{eqnarray*}
\Big(\int_M|W+\frac{\sqrt{2}}{\sqrt{n}(n-2)}Z \KN g|^2 \Big)^{\frac{2}{n}}&&+\sqrt{\frac{(n-4)^2(n-1)}{8(n-2)}}\lambda vol(M)^{\frac{2}{n}}\\
&&\leq \sqrt{\frac{n-2}{32(n-1)}}Y(M,[g]),
\end{eqnarray*} 
where $Z=Ric-\frac{R}{n}g$ is the trace-less Ricci tensor, $W$ is the Weyl tensor and $Y(M,[g])$ is the Yamabe invariant associated to $(M,g)$, then $M$ is isometric to a quotient of the round sphere $\mathbb{S}^n$.
\end{theorem}
 In 2019, Mondal and Shaikh \cite{CA19} proved the isometry theorem for gradient $\rho$-Einstein soliton in case of conformal vector field. In particular, they proved the following result:
\begin{theorem}\cite{CA19}
Let $(M,g,f,\rho)$ be a compact gradient $\rho$-Einstein soliton. If $\nabla f$ is a non-trivial conformal vector field, then $M$ is isometric to the Euclidean sphere $\mathbb{S}^n$.
\end{theorem}
Dwivedi \cite{DW18} proved an isometry theorem for gradient Ricci-Bourguignon soliton. 
\begin{theorem}\label{th1}\cite{DW18}
A non-trivial compact gradient Ricci-Bourguignon soliton is isometric to an Euclidean sphere if any one of the following holds\\
(1) $M$ has constant scalar curvature.\\
(2) $\int_Mg(\nabla R,\nabla f)\leq 0$.\\
(3) $M$ is a homogeneous manifold.
\end{theorem}
We note that Catino et. al. \cite{CMM15} proved many results for gradient $\rho$-Einstein soliton in non-compact manifold.  
\begin{theorem}
Let $(M,g,f,\rho)$ be a complete non-compact gradient shrinking $\rho$-Einstein soliton with $0<\rho<1/2(n-1)$ bounded curvature, non-negative radial sectional curvature, and non-negative Ricci curvature. Then the scalar curvature is constant.
\end{theorem}
In this paper, we have showed that a non-trivial compact gradient $\rho$-Einstein soliton is isometric to an Euclidean sphere. The main results of this paper are as follows:
\begin{theorem}\label{th2}
A nontrivial compact gradient $\rho$-Einstein soliton has constant scalar curvature and therefore $M$ is isometric to an Euclidean sphere.
\end{theorem}
We have also showed that in a non-compact  gradient $\rho$-Einstein soliton satisfying some conditions the scalar curvature vanishes.
\begin{theorem}\label{th3}
Suppose $(M,g,f,\rho)$ is a non-compact gradient non-expanding $\rho$-Einstein soliton with non-negative scalar curvature. If $\rho> 1/n$ and the $\rho$-Einstein potential function satisfies
 \begin{equation}
 \int_{M-B(p,r)}d(x,p)^{-2}f<\infty,
 \end{equation}
 then the scalar curvature vanishes in $M$.
\end{theorem}
\section{Proof of the results}
\begin{proof}[\textbf{Proof of the Theorem \ref{th2}}]
Since the gradient $\rho$-Einstein soliton is non-trivial, it follows that $\rho\neq 1/n$, see \cite{CM16}. Taking the trace of (\ref{eq1}) we get
\begin{equation}\label{eq2}
R+\Delta f=\lambda n+\rho Rn.
\end{equation}
From the commutative equation, we obtain
\begin{equation}\label{eq3}
\Delta\nabla_if=\nabla_i\Delta f+R_{ij}\nabla_jf.
\end{equation}
By using contracted second Bianchi identity, we have
\begin{eqnarray*}
\Delta\nabla_if=\nabla_j\nabla_j\nabla_i f&=&\nabla_j(\lambda g_{ij}+\rho Rg_{ij}-R_{ij})\\
&=& \nabla_i(\rho R-\frac{1}{2}R).
\end{eqnarray*}
and
$$\nabla_i\Delta f=\nabla_i(\lambda n+\rho Rn-R)=\nabla_i(\rho Rn-R).$$
Therefore, (\ref{eq3}) yields
\begin{equation}
(n-1)\rho \nabla_iR-\frac{1}{2}\nabla_iR+R_{ij}\nabla_jf=0,
\end{equation}
Taking covariant derivative $\nabla_l$, we get
$$(n-1)\rho \nabla_l\nabla_iR-\frac{1}{2}\nabla_l\nabla_iR+\nabla_lR_{ij}\nabla_jf+R_{ij}\nabla_l\nabla_jf=0.$$
Taking trace in both sides, we obtain 
\begin{equation}
((n-1)\rho-\frac{1}{2})\Delta R+\frac{1}{2}g(\nabla R,\nabla f)+R(\lambda n+\rho Rn-R)=0.
\end{equation}
Now integrating using divergence theorem we get
\begin{eqnarray*}
\int_M R(\lambda n+\rho Rn-R)&=&-\int_M ((n-1)\rho-\frac{1}{2})\Delta R-\frac{1}{2}\int_M g(\nabla R,\nabla f)\\
&=& \frac{1}{2}\int_M R\Delta f=\frac{1}{2}\int_M R(\lambda n+\rho Rn-R).
\end{eqnarray*}
The above equation is true only if 
\begin{equation}
\int_M R(\lambda n+\rho Rn-R)=0,
\end{equation}
which implies 
\begin{equation}\label{eq4}
\int_M R\Big(R+\frac{\lambda n}{n\rho-1}\Big)=0,
\end{equation}
Again integrating (\ref{eq2}), we obtain
\begin{equation}\label{eq5}
\int_M \Big(R+\frac{\lambda n}{n\rho-1}\Big)=0.
\end{equation}
Therefore, (\ref{eq4}) and (\ref{eq5}) together imply that
$$\int_M \Big(R+\frac{\lambda n}{n\rho-1}\Big)^2=0.$$
Hence, $R=\lambda n/(1-\rho n)$. Then from Theorem \ref{th1} we can conclude our result. 
\end{proof}

\begin{proof}[\textbf{Proof of the Theorem \ref{th3}}]
From (\ref{eq2}) we get
$$(n\rho-1)R=\Delta f-\lambda n.$$
Since $\lambda\geq 0$, the above equation implies that
\begin{equation}
(n\rho-1)R\leq  \Delta f.
\end{equation}
Now, we consider the cut-off function, introduced in \cite{CC96}, $\varphi_r\in C^2_0(B(p,2r))$ for $r>0$ such that
\[ \begin{cases} 
	  0\leq \varphi_r\leq 1 &\text{ in }B(p,2r)\\
      \varphi_r=1  & \text{ in }B(p,r) \\
      |\nabla \varphi_r|^2\leq\frac{C}{r^2}& \text{ in }B(p,2r) \\
      \Delta \varphi_r\leq \frac{C}{r^2} &  \text{ in }B(p,2r),
   \end{cases}
\]
where $C>0$ is a constant. Then for $r\rightarrow\infty$, we have $\Delta \varphi^2_r\rightarrow 0$ as $\Delta \varphi^2_r\leq \frac{C}{r^2}$. Then we calculate
\begin{eqnarray}
(n\rho-1)\int_MR\varphi^2_r\leq \int_M \varphi^2_r\Delta f&=&\int_{B(p,2r)-B(p,r)}f\Delta \varphi_r^2\\
&\leq & \int_{B(p,2r)-B(p,r)}f\frac{C}{r^2}\rightarrow 0,
\end{eqnarray}
 as $r\rightarrow \infty$. Hence, we obtain
 \begin{equation}
 (n\rho-1)\lim_{r\rightarrow \infty}\int_{B(p,r)}R\leq 0.
 \end{equation}
Since $\rho>1/n$, it follows that 
$$\lim_{r\rightarrow \infty}\int_{B(p,r)}R\leq 0.$$
But $R$ is non-negative everywhere in $M$. Therefore, $R\equiv 0$ in $M$. 
\end{proof}
\section{acknowledgment}
 The third author gratefully acknowledges to the
 CSIR(File No.:09/025(0282)/2019-EMR-I), Govt. of India for financial assistance.

\end{document}